\documentclass[11pt]{article}

\usepackage[T1]{fontenc}
\usepackage[utf8]{inputenc}
\usepackage{lmodern}
\usepackage{microtype}
\usepackage[margin=1in]{geometry}
\usepackage{amsmath,amssymb,amsthm,mathtools}
\usepackage{enumitem}
\usepackage[ pdfusetitle,
             colorlinks=true,
             linkcolor=blue,
             citecolor=blue,
             urlcolor=blue ]{hyperref}

\setlist{topsep=3pt, itemsep=2pt, parsep=0pt}
\newtheorem{theorem}{Theorem}
\newtheorem{corollary}[theorem]{Corollary}
\newtheorem{lemma}[theorem]{Lemma}
\newtheorem{claim}[theorem]{Claim}
\usepackage{cleveref}
\crefname{claim}{Claim}{Claims}
\usepackage{natbib}

\newcommand{\R}{\mathbb{R}}
\newcommand{\C}{\mathbb{C}}

\newcommand{\Tr}{\operatorname{tr}}
\newcommand{\tr}{\Tr}
\newcommand{\T}{\mathsf{T}}
\newcommand{\Sym}{\mathcal{S}}
\newcommand{\Herm}{\mathcal{H}}
\newcommand{\rmre}{\mathrm{re}}
\newcommand{\rmim}{\mathrm{im}}

\newcommand{\ignore}[1]{}
\newcommand{\mnote}[1]{[\texttt{{\color{blue} #1}]}}

\title{A Non-Commutative Binomial Theorem, Extensions, and Applications}
\author{Ravi Kumar \\
Google Research \\
Mountain View, CA \\
\texttt{ravi.k53@gmail.com}
\and
Raghu Meka \\
Google DeepMind \\
Bengaluru, India \\
\texttt{raghuvardhan@gmail.com}
}
\date{}

\begin{document}
\maketitle

\begin{abstract}
This paper presents a matrix trace inequality analogous to the binomial theorem for reals: for Hermitian matrices $A$ and $B$, and a nonnegative even integer $k$ (and any nonnegative integer $k$ when $A, B$ are positive semi-definite),
\[
\Tr\bigl((A + B)^k\bigr) \le \sum_{\ell=0}^k \binom{k}{\ell}\Tr\bigl(A^\ell B^{k-\ell}\bigr).
\] 
The above natural inequality implies the classical Golden--Thompson inequality $\Tr(e^{A + B}) \le \Tr(e^A e^B)$.  It also proves an extension, which can be used to show other new inequalities as well as give short proofs of inequalities such as the Schatten-$p$ two-point inequality for $p \ge 3$ and Hanner's inequality for $p \ge 3$.

The main proof was discovered by Gemini Carbon within an agentic system.
\end{abstract}

\section{Introduction}
The main result of this paper is the following trace inequality that can be seen as a non-commutative version of the usual binomial theorem over scalars:
\begin{theorem}[Matrix binomial theorem]
\label{thm:binomial}
Let $A, B$ be $d\times d$ Hermitian matrices. Then, for every even integer $k$, 
\begin{equation}
\label{eq:binomial}
\Tr\bigl((A + B)^k\bigr) \le \sum_{\ell=0}^k \binom{k}{\ell}\Tr\bigl(A^\ell B^{k-\ell}\bigr).
\end{equation}
Moreover, if $A,B$ are positive semi-definite, then the above holds for all integers $k \geq 0$. 
%with equality for $k \in \{0, 1,2,3\}$.
\end{theorem}

\ignore{
\begin{theorem}[Matrix Binomial Theorem]
\label{thm:binomial}
For positive semi-definite matrices $A, B \in \C^{d \times d}$ and any even integer $k \geq 2$, 
\begin{equation}
\label{eq:binomial}
\Tr\bigl((A + B)^k\bigr) \le \sum_{\ell=0}^k \binom{k}{\ell}\Tr\bigl(A^\ell B^{k-\ell}\bigr).
\end{equation}
\end{theorem}}

Note that standard properties of trace imply that the inequality in \Cref{thm:binomial} holds with equality for $k \in \{0,1,2,3\}$. Further, for general Hermitian matrices, the inequality can fail for odd $k \geq 5$. Throughout, we use the following notation: for an integer $d \geq 1$, 
\begin{itemize}[nosep]
\item For $x \in \C^d$, 
let $x^*$ be its conjugate transpose. 
\item Let $I_d \in \R^{d \times d}$ denote the identity matrix.
\item For $A \in \C^{d \times d}$, let 
$\Tr(A) := \sum_{i=1}^d A_{ii}$,
$A^*$ be its conjugate transpose, and $A^0 := I_d$.
\item Let $\Herm^d := \{A \in \C^{d \times d} \mid A^* = A\}$ and $\Herm_+^d := \{A \in \Herm^d \mid x^* A x \ge 0\text{ for all }x \in \C^d\}$ 
be the sets of Hermitian and positive semi-definite (PSD)  matrices. 
\end{itemize}
As a motivation, note that \Cref{thm:binomial} immediately implies the classical Golden--Thompson inequality: \begin{corollary}[Golden--Thompson inequality~\citep{golden, thompson}]
\label{cor:gt}
For all $A, B \in \Herm^d$,
\[
\Tr\bigl(e^{A + B}\bigr) \le \Tr\bigl(e^A e^B\bigr).
\]
\end{corollary}
\begin{proof}
Let $c \ge 0$ be large enough such that $\tilde{A} := A + c I_d \in \Herm_+^d$ and $\tilde{B} := B + c I_d \in \Herm_+^d$.  Now, 
\begin{align*}
e^{2c}\Tr\bigl(e^{A + B}\bigr) = \Tr\bigl(e^{\tilde{A} + \tilde{B}}\bigr) &= \sum_{k=0}^\infty \frac{1}{k!}\Tr\bigl((\tilde{A} + \tilde{B})^k\bigr) 
\\
\mbox{(using Theorem~\ref{thm:binomial})
} &\le \sum_{k=0}^\infty \sum_{\ell=0}^k \frac{1}{\ell!(k - \ell)!}\Tr\bigl(\tilde{A}^\ell \tilde{B}^{k-\ell}\bigr) = \Tr\bigl(e^{\tilde{A}} e^{\tilde{B}}\bigr) = e^{2c}\Tr\bigl(e^A e^B\bigr),
\end{align*}
and dividing by $e^{2c} > 0$ completes the proof. 
\end{proof}

\paragraph{Related Work}
\Cref{thm:binomial}, while quite natural, surprisingly does not seem to be known or explicitly\footnote{The second author was led to this exact statement in 2013 while pursuing a now defunct application; several experts at the time could not point to any concrete references either.} asked in the literature to the best of our knowledge. However, several related forms and strengthenings have been studied, as discussed below.

Expanding $(A + B)^k$ gives the sum of all $2^k$ words of length $k$ in $A$ and $B$. Let $S_{\ell,k-\ell}(A, B)$ denote the sum of the $\binom{k}{\ell}$ words containing $\ell$ copies of $A$ and $k - \ell$ copies of $B$, so that $\tr((A + B)^k) = \sum_{\ell=0}^k \tr(S_{\ell,k-\ell}(A, B))$. For $A, B \in \mathcal{H}_+^d$, the \emph{Bessis--Moussa--Villani (BMV) conjecture}~\citep{BMV1975}, in the equivalent polynomial formulation of~\cite{LiebSeiringer2004}, asserts that $\tr(S_{\ell,k-\ell}(A, B)) \ge 0$ for every $\ell$; this was proved by~\cite{Stahl2013}. H\"agele's \emph{refined BMV conjecture}~\citep{HageleOQP40} asked whether each slice also satisfies the upper bound $\tr(S_{\ell,k-\ell}(A, B)) \le \binom{k}{\ell}\tr(A^\ell B^{k-\ell})$, which summed over $\ell$ would immediately yield \Cref{thm:binomial} on $\mathcal{H}_+^d$. While alternating-product inequalities of~\cite{ando-hiai-okubo} imply this slice bound on $\mathcal{H}_+^d$ for all $k \le 8$,~\cite{plevnik} disproved the individual word bound $\tr(w(A, B)) \le \tr(A^\ell B^{k-\ell})$ for $3 \times 3$ matrices in $\mathcal{H}_+^3$ at $(k, \ell) = (10, 5)$, and~\cite{ChaLee2026} (see also~\citep{Dinh2026}) recently disproved the slice bound itself on $\mathcal{H}_+^3$ at $\ell = k - \ell = 5$. Thus cross-slice cancellation is essential, and \Cref{thm:binomial} shows that once one sums over all $\ell \in \{0, \dots, k\}$, the upper bound holds for every $k$.

A second connection comes from expanding both sides of the Golden--Thompson inequality $\tr(e^{t(A+B)}) \le \tr(e^{tA} e^{tB})$~\citep{golden,thompson} as power series in $t$: the coefficients of $t^k / k!$ on the two sides are $\tr((A + B)^k)$ and $\sum_{\ell=0}^k \binom{k}{\ell}\tr(A^\ell B^{k-\ell})$, so \Cref{thm:binomial} strengthens Golden--Thompson coefficient by coefficient on $\mathcal{H}_+^d$ (and for all even degrees on $\mathcal{H}^d$). Equality holds for $k \in \{0, 1, 2, 3\}$, while for $k = 4$ the gap is given by the classical identity $2\tr(A^2 B^2 - ABAB) = -\tr([A, B]^2) = \|[A, B]\|_F^2 \ge 0$ (see, e.g.,~\cite{Tao2010GT}), where $\|X\|_F := \tr(X^* X)^{1/2}$. 

Recently,~\cite{zhang} studied eigenvalue comparisons for these Taylor coefficients on $\mathcal{H}^d$. Setting $Q_k := \sum_{\ell=0}^k \binom{k}{\ell} A^\ell B^{k-\ell}$ and $\operatorname{Re} Q_k := \frac{1}{2}(Q_k + Q_k^*)$, that paper proved that for $k \in \{3, 4\}$ the sum of the $m$ largest eigenvalues of $(A + B)^k$ is at most that of $\operatorname{Re} Q_k$ for every $1 \le m \le d$ (weak majorization, written $\lambda((A + B)^k) \prec_w \lambda(\operatorname{Re} Q_k)$), and conjectured that this holds for all $A, B \in \mathcal{H}^d$ and all $k \ge 5$~\cite[Conjecture~5.5]{zhang}. Since $\tr(Q_k^*) = \tr(Q_k)$, the $m = d$ case of that conjecture is the trace inequality $\tr((A + B)^k) \le \tr(Q_k)$. Moreover, trace cyclicity and Vandermonde's identity $\sum_{r+s=\ell}\binom{k}{r}\binom{k}{s} = \binom{2k}{\ell}$ give $\|Q_k\|_F^2 = \tr(Q_k^* Q_k) = \tr(Q_{2k})$, so the even-exponent case of this trace inequality is equivalent to the Frobenius-norm bound $\|(A + B)^k\|_F \le \|Q_k\|_F$ for all $k \geq 0$. \Cref{thm:binomial} establishes $\tr((A + B)^k) \le \tr(Q_k)$ for all even $k \geq 0$ on $\mathcal{H}^d$ (hence $\|(A + B)^k\|_F \le \|Q_k\|_F$ for all $k \geq 0$) and for all $k \geq 0$ on $\mathcal{H}_+^d$, whereas for odd $k \ge 5$ and general $A, B \in \mathcal{H}^d$ it---and hence~\cite[Conjecture~5.5]{zhang} as stated---already fails for $2 \times 2$ matrices.\footnote{For instance, $A = -\operatorname{diag}(1, 0)$ and $B = -\bigl(\begin{smallmatrix}1 & 1 \\ 1 & 1\end{smallmatrix}\bigr)$ give $\tr((A + B)^5) = -123 > -138 = \tr(Q_5)$.}

Finally, the proof method here yields a more general trace inequality (see \Cref{thm:general}): for $A, B \in \mathcal{H}^d$ with eigendecompositions $A = \sum_{i=1}^d \lambda_i u_i u_i^*$ and $B = \sum_{j=1}^d \mu_j v_j v_j^*$ and eigenbasis \emph{overlaps} $W_{ij} := |u_i^* v_j|^2$, we have
\[
\tr\bigl(f(A + B)\bigr) \le \sum_{i,j=1}^d W_{ij} f(\lambda_i + \mu_j)
\]
whenever $f''$ is convex on the relevant spectral interval---in particular, for every $f \in C^4(\mathbb{R})$ with $f^{(4)} \ge 0$ (taking $f(x) = x^k$ for $k$ even recovers \Cref{thm:binomial}). Recent work of~\cite{Heinavaara2026} established an asymptotic counterpart of this phenomenon in free probability. Given compactly supported probability measures $\mu$ and $\nu$ on $\mathbb{R}$, let $\mu * \nu$ denote their classical convolution (the law of $X + Y$ for independent $X \sim \mu$, $Y \sim \nu$) and let $\mu \boxplus \nu$ denote their \emph{free additive convolution}---which arises as the $d \to \infty$ limiting eigenvalue distribution of $A + U B U^*$ when $A, B \in \mathcal{H}^d$ have empirical eigenvalue distributions approaching $\mu, \nu$ and $U$ is a Haar-random $d \times d$ unitary matrix. Hein\"avaara proved that $\int f\,d(\mu \boxplus \nu) \le \int f\,d(\mu * \nu)$ holds for all such $\mu, \nu$ if and only if $f^{(4)} \ge 0$ (in particular ordering all even moments).

The necessity of $f^{(4)} \ge 0$ follows from a short perturbative computation (comparing $\mu * \varepsilon\nu$ with $\mu \boxplus \varepsilon\nu$ as $\varepsilon \to 0$; see~\cite[Section~3]{Heinavaara2026}).
%The necessity of $f^{(4)} \ge 0$ can be shown by a direct example (taking $\mu = \nu$ to be a two-point measure of small width as in~\cite[Section~3]{Heinavaara2026}). 
The substance of the theorem is the inequality itself, whose proof in~\citep{Heinavaara2026} constructs a positive measure on $\mathbb{R}^2$ using free cumulants and identities for Gegenbauer polynomials. Because this statement averages over the relative position of the eigenbases of $A$ and $B$, it does not imply \Cref{thm:binomial} for fixed matrices. Its proof also does not seem to transfer directly.%: it yields a separate comparison for each slice, and for fixed matrices such slice-wise bounds can fail~\citep{ChaLee2026}.

%Because this statement averages out the relative eigenbases of $A$ and $B$, it does not imply \Cref{thm:binomial} for fixed matrices; moreover, it yields slice-by-slice polynomial inequalities, which do not extend to fixed matrices in view of~\citep{ChaLee2026}.

In the other direction, applying \Cref{thm:general} to $A$ and $U B U^*$, averaging over a Haar-random unitary $U$ (so that $\mathbb{E}_U[|u_i^* U v_j|^2] = \frac{1}{d}$), and letting $d \to \infty$ recovers Hein\"avaara's comparison inequality, and hence the sufficiency direction of~\cite[Theorem~1.3]{Heinavaara2026}, by a different route (see \Cref{sec:free}). It does not recover the positive measure on $\mathbb{R}^2$ or the term-by-term comparisons from~\citep{Heinavaara2026}.

\paragraph{AI Acknowledgements} The main proof was found by Gemini Carbon within an agentic system. The original proof was quite long and involved advanced tools. The authors simplified the proof in many further interactions, bringing it to its present self-contained short form while also identifying sources of some of the ideas in existing works. The generalization to other functions (\Cref{thm:general}) was suggested by the authors. The applications were discovered in further interactions with the agent and surveying the known literature.  

\section{Main Proof}
\label{sec:main}

For $s \in \R$, let $s_+ := \max(0, s)$
and $s_+^p := (s_+)^p$ for
$p > 0$.
For $x \in \C^d$, 
let $\|x\|_2 := (x^* x)^{1/2}$. 
For any orthonormal basis $(x_1, \dots, x_d)$ of $\C^d$ and any $x \in \C^d$, $x = \sum_{i=1}^d (x_i^* x) x_i$, so $\sum_{i=1}^d x_i x_i^* = I_d$ and $\|x\|_2^2 = \sum_{i=1}^d |x_i^* x|^2$. 

Let $A \in \Herm^d$ have  eigendecomposition $A = \sum_{i=1}^d \lambda_i u_i u_i^*$, where $\lambda_i \in \R$ and the $u_i$'s are orthonormal.  For any real-valued function $g$ defined
on the $\lambda_i$'s, consider the matrix $g(A) := \sum_{i=1}^d g(\lambda_i)u_i u_i^* \in \Herm^d$; it is easy to see that the matrix $g(A)$ is well defined.  In particular, $A_+ := \sum_{i=1}^d (\lambda_i)_+ u_i u_i^* \in \Herm_+^d$; note also that $\|A_+ x\|_2^2 = \sum_{i=1}^d (\lambda_i)_+^2 |u_i^* x|^2 \le \sum_{i=1}^d \lambda_i^2 |u_i^* x|^2 = \|A x\|_2^2$ for all $x \in \C^d$.
Moreover, $\Tr(g(A)) = \sum_{i=1}^d g(\lambda_i)$.

For the rest of this section, fix matrices $A, B \in \Herm^d$, let $A = \sum_{i=1}^d \lambda_i u_i u_i^*$ and $B = \sum_{j=1}^d \mu_j v_j v_j^*$ be their respective eigendecompositions, and let $W_{ij} := |u_i^* v_j|^2$ be the \emph{eigenbasis overlaps}. In case of ambiguity, the notation will be called out explicitly. 

First, note that the RHS of \Cref{thm:binomial} can be rewritten as follows by expanding $A, B$ in their eigenbases:

\begin{lemma}
\label{lem:overlap}
(i)~$\sum_{\ell=0}^k \binom{k}{\ell}\tr(A^\ell B^{k-\ell}) = \sum_{i,j=1}^d W_{ij}(\lambda_i + \mu_j)^k$ for every integer $k \geq 0$.
(ii)~$\tr((A + B)^k) = \sum_{i,j=1}^d W_{ij}(\lambda_i + \mu_j)^k$ for $k \in \{0, 1, 2, 3\}$.
\end{lemma}
\begin{proof}
For (i), we have $A^\ell = \sum_{i=1}^d \lambda_i^\ell u_i u_i^*$ and $B^{k-\ell} = \sum_{j=1}^d \mu_j^{k-\ell} v_j v_j^*$. Using $\tr(u_i u_i^* v_j v_j^*) = \tr(v_j^* u_i u_i^* v_j) = |u_i^* v_j|^2 = W_{ij}$ and the binomial theorem for reals, we obtain
\[
\sum_{\ell=0}^k \binom{k}{\ell}
\tr(A^\ell B^{k-\ell})
= \sum_{i,j=1}^d W_{ij}\sum_{\ell=0}^k \binom{k}{\ell}\lambda_i^\ell \mu_j^{k-\ell}
= \sum_{i,j=1}^d W_{ij}(\lambda_i + \mu_j)^k.
\]

For (ii), expanding $(A + B)^k$ for $k \in \{0, 1, 2, 3\}$ and using $\tr(XY) = \tr(YX)$ (in particular, $\tr(BA) = \tr(AB)$, $\tr(ABA) = \tr(A^2 B)$, and $\tr(BAB) = \tr(A B^2)$) gives $\tr((A + B)^0) = \tr(I_d) = \tr(A^0 B^0)$, 
$\tr(A + B) = \tr(B) + \tr(A)$, $\tr((A + B)^2) = \tr(B^2) + 2\tr(AB) + \tr(A^2)$, and $\tr((A + B)^3) = \tr(B^3) + 3\tr(A B^2) + 3\tr(A^2 B) + \tr(A^3)$. Thus $\tr((A + B)^k) = \sum_{\ell=0}^k \binom{k}{\ell}\tr(A^\ell B^{k-\ell})$ for $k \in \{0, 1, 2, 3\}$.  Finally,
(ii) follows from (i). 
\end{proof}

Thus, if we take $f: \R \to \R$ to be defined as $f(x) = x^k$, then \Cref{thm:binomial} is equivalent to proving
\[
\tr(f(A+B)) \leq \sum_{i,j=1}^d W_{ij} f(\lambda_i + \mu_j).
\]

The above version stresses that it is important to exploit how the ``overlap'' coefficients $W_{ij}$ and the eigenvalues of $A+B$ relate to each other. The key technical ingredient is that instead of targeting the above inequality for the function $f(x) = x^k$ directly, it is easier to first prove it for the map $x \mapsto x_+^3$. This strategy is reminiscent of the approach of~\cite{AizenmanLieb1978}, who use a similar integral identity to control higher powers of eigenvalues of Schr\"odinger operators. 

\begin{lemma}[Main Technical Lemma]\label{lm:maintech}
For $A, B \in \Herm^d$ with eigendecompositions $A = \sum_{i=1}^d \lambda_i u_i u_i^*$ and $B = \sum_{j=1}^d \mu_j v_j v_j^*$, $W_{ij} = |u_i^* v_j|^2$,
\begin{equation}\label{eq:peano}
\tr((A+B)^3_+) \leq \sum_{i,j=1}^d W_{ij} (\lambda_i + \mu_j)^3_+.    
\end{equation}
\end{lemma}

The above lemma implies the main binomial theorem by a careful integration by parts. 

\begin{proof}[Proof of \Cref{thm:binomial} from \Cref{lm:maintech}]
By Lemma~\ref{lem:overlap}(ii), equality holds in~\eqref{eq:binomial} for $k \in \{0, 1, 2, 3\}$.  Let $\sigma_1, \dots, \sigma_d$ be the eigenvalues of $A + B$, so $\tr((A + B)^k) = \sum_{i=1}^d \sigma_i^k$ and $\tr((A + B - tI_d)_+^3) = \sum_{i=1}^d (\sigma_i - t)_+^3$.  
Let $T := \max \big( \{|\lambda_i + \mu_j|\}_{i,j=1}^d \cup \{|\sigma_i|\}_{i=1}^d \big) \ge 0$.  
For every integer $k \ge 4$ and $x \in [0, T]$, integrating by parts three times gives
\begin{equation}
\label{eq:beta}
\int_0^T (x - t)_+^3 t^{k-4}\, dt
= \int_0^x (x - t)^3 t^{k-4}\, dt
= \frac{6}{k(k - 1)(k - 2)(k - 3)}\, x^k
= \frac{1}{4\binom{k}{4}}\, x^k.
\end{equation}
Define the function
\[
D(t; A, B) := \sum_{i,j=1}^d W_{ij}(\lambda_i + \mu_j - t)_+^3 - \sum_{i=1}^d (\sigma_i - t)_+^3,
\]
which is continuous in $t$ and by \Cref{lm:maintech} applied to $A, B - t I_d$, satisfies $D(t; A, B) \ge 0$ for all $t \in \R$.

When $A, B \in \Herm_+^d$, we also have $A + B \in \Herm_+^d$, so $\lambda_i + \mu_j \in [0, T]$ and $\sigma_i \in [0, T]$ for all $i, j$. Applying~\eqref{eq:beta} to $x = \lambda_i + \mu_j$ and to $x = \sigma_i$, and using Lemma~\ref{lem:overlap}(i), yields for every $k \ge 4$,
\[
\sum_{\ell=0}^k \binom{k}{\ell}\tr(A^\ell B^{k-\ell}) - \tr\bigl((A + B)^k\bigr)
= \sum_{i,j=1}^d W_{ij}(\lambda_i + \mu_j)^k - \sum_{i=1}^d \sigma_i^k
= 4\binom{k}{4}\int_0^T D(t; A, B)\, t^{k-4}\, dt \ge 0.
\]

When $A, B \in \Herm^d$ and $k \ge 4$ is even, any $y \in [-T, T]$ satisfies $|y| \in [0, T]$, $y^{k} = |y|^{k}$, and $(|y| - t)_+^3 = (y - t)_+^3 + (-y - t)_+^3$ for all $t \ge 0$, so applying~\eqref{eq:beta} to $x = |y|$ yields
\[
y^{k} = 4\binom{k}{4}\int_0^T \bigl((y - t)_+^3 + (-y - t)_+^3\bigr) t^{k-4}\, dt.
\]
Since $-A$, $-B$, and $-(A + B)$ have eigenvalues $-\lambda_i$, $-\mu_j$, and $-\sigma_i$ with the same overlap weights $W_{ij}$, applying this identity to $y = \lambda_i + \mu_j$ and to $y = \sigma_i$ yields
\[
\sum_{\ell=0}^{k}\binom{k}{\ell}\tr(A^\ell B^{k-\ell}) - \tr\bigl((A + B)^{k}\bigr)
= 4\binom{k}{4}\int_0^T \bigl( D(t; A, B) + D(t; -A, -B)\bigr) t^{k-4}\, dt \ge 0.\qedhere
\]
\end{proof}

\subsection{Proof of \texorpdfstring{\Cref{lm:maintech}}{Lemma~\ref{lm:maintech}}}
The proof of \Cref{lm:maintech} relies on the following subgradient inequality: 

\begin{lemma}
\label{lem:subgrad}
For all $X, Y \in \Herm^d$,
\[
\frac{1}{3}\Tr(Y_+^3) - \frac{1}{3}\Tr(X_+^3) \le \Tr\bigl(Y_+^2 (Y - X)\bigr).
\]
\end{lemma}
This is Klein's inequality $\Tr(f(Y)) - \Tr(f(X)) \le \Tr\bigl(f'(Y)(Y - X)\bigr)$~\citep{Klein1931} (see, e.g.,~\cite{Carlen2010}) for the convex function $f(x) = \frac{1}{3}x_+^3$; a
self-contained proof for this special case is included below for completeness.
\begin{proof}
Consider the eigendecompositions $X = \sum_{i=1}^d \lambda_i u_i u_i^*$ and $Y = \sum_{j=1}^d \mu_j v_j v_j^*$. For all $a, b \ge 0$, we have the inequality 
\[
\frac{2}{3}a^3 + \frac{1}{3}b^3 - a^2 b = \frac{1}{3}(a - b)^2(2a + b) \ge 0.
\]

Using the above with $a = (\mu_j)_+$ and $b = (\lambda_i)_+$, and $\sum_{i=1}^d |u_i^* v_j|^2 = \sum_{j=1}^d |u_i^* v_j|^2 = 1$, we obtain
\[
\tr(Y_+^2 X) = \sum_{i,j=1}^d |u_i^* v_j|^2 (\mu_j)_+^2 \lambda_i
\le \sum_{i,j=1}^d |u_i^* v_j|^2 \left(\frac{2}{3}(\mu_j)_+^3 + \frac{1}{3}(\lambda_i)_+^3\right)
= \frac{2}{3}\tr(Y_+^3) + \frac{1}{3}\tr(X_+^3).
\]
Subtracting this from $\tr(Y_+^2 Y) = \tr(Y_+^3)$ yields the lemma.
\end{proof}
Towards the proof of \Cref{lm:maintech}, define $y_{\max} := \max_{1 \le j \le d} \mu_j$, and set $x_i := \lambda_i + y_{\max}$ and $y_j := \mu_j - y_{\max}$; by definition, $x_i + y_j = \lambda_i + \mu_j$ and $y_j \le 0$. Define $C := A + B = \sum_{i=1}^d x_i u_i u_i^* + \sum_{j=1}^d y_j v_j v_j^*$. Let $s_{\max} := \max(0, \max_{1 \le i \le d} x_i) \ge 0$, and for $s \in [0, s_{\max}]$ define $x_i(s) := \min(\max(x_i, 0), s)$, $y_j(s) := \max(y_j, -s)$, and the functions
\begin{align*}
C(s) &:= \sum_{i=1}^d x_i(s) u_i u_i^* + \sum_{j=1}^d y_j(s) v_j v_j^*, \\
L(s) &:= \frac{1}{3}\tr\bigl(C(s)_+^3\bigr), \\
R(s) &:= \frac{1}{3}\sum_{i,j=1}^d W_{ij}\bigl(x_i(s) + y_j(s)\bigr)_+^3.
\end{align*}
Note that $L(0) = R(0) = 0$. 

\paragraph{Proof idea.} We reach $A + B$ through the family $C(s)$, in which the coefficients $x_i$ are capped at $s$ and the coefficients $y_j$ are floored at $-s$. Both $L$ and $R$ vanish at $s = 0$. \Cref{clm:Lstep,clm:Rstep} show that when $s$ increases by a small step $h$, $L$ grows by at most as much as $R$, up to an $O(h^2)$ error; summing over steps gives $L(s_{\max}) \le R(s_{\max})$. \Cref{clm:endpoint} relates the endpoint $s = s_{\max}$ to $A + B$.
\begin{claim}
\label{clm:endpoint}
$\frac{1}{3}\tr(C_+^3) \le L(s_{\max})$ and $R(s_{\max}) = \frac{1}{3}\sum_{i,j=1}^d W_{ij}(\lambda_i + \mu_j)_+^3$.
\end{claim}

\begin{proof}
At $s = s_{\max}$, we have $x_i \le x_i(s_{\max}) = \max(x_i, 0)$ and $y_j \le y_j(s_{\max}) = \max(y_j, -s_{\max})$. Applying Lemma~\ref{lem:subgrad} to $X = C(s_{\max}), Y = C$ and using $\tr(C_+^2 z z^*) = \|C_+ z\|_2^2$ yields
\begin{align*}
\frac{1}{3}\tr(C_+^3) - L(s_{\max})
&\le \tr\bigl(C_+^2(C - C(s_{\max}))\bigr) \\
&= \tr\left(C_+^2 \left(\sum_{i=1}^d (x_i - x_i(s_{\max})) u_i u_i^*\right)\right) + \tr\left(C_+^2 \left(\sum_{j=1}^d (y_j - y_j(s_{\max})) v_j v_j^*\right)\right) \\
&= \sum_{i=1}^d (x_i - x_i(s_{\max}))\tr(C_+^2 u_i u_i^*) + \sum_{j=1}^d (y_j - y_j(s_{\max}))\tr(C_+^2 v_j v_j^*) \\
&= \sum_{i=1}^d (x_i - x_i(s_{\max}))\|C_+ u_i\|_2^2 + \sum_{j=1}^d (y_j - y_j(s_{\max}))\|C_+ v_j\|_2^2 \\
&\le 0.
\end{align*}
Next, for each $(i, j)$: if $x_i < 0$, then $x_i(s_{\max}) = 0$ and $y_j(s_{\max}) \le 0$, so $(x_i(s_{\max}) + y_j(s_{\max}))_+ = 0 = (x_i + y_j)_+$ since $y_j \le 0$. If $y_j < -s_{\max}$, then $y_j(s_{\max}) = -s_{\max}$, so $(x_i(s_{\max}) + y_j(s_{\max}))_+ = 0 = (x_i + y_j)_+$ since $x_i \le x_i(s_{\max}) \le s_{\max}$. Otherwise, $x_i \ge 0$ and $y_j \ge -s_{\max}$, so $x_i(s_{\max}) = x_i$ and $y_j(s_{\max}) = y_j$. Thus, for all $(i, j)$, $(x_i(s_{\max}) + y_j(s_{\max}))_+ = (x_i + y_j)_+ = (\lambda_i + \mu_j)_+$, giving $R(s_{\max}) = \frac{1}{3}\sum_{i,j=1}^d W_{ij}(\lambda_i + \mu_j)_+^3$.
\end{proof}

Now suppose $s_{\max} > 0$, fix an integer $N \ge 1$, and order the distinct elements of $\{\frac{n}{N}s_{\max}\}_{n=0}^N \cup \{\max(x_i, 0)\}_{i=1}^d \cup \{\min(-y_j, s_{\max})\}_{j=1}^d \subseteq [0, s_{\max}]$ as $0 = \tau_0 < \tau_1 < \dots < \tau_K = s_{\max}$. Fix $n \in \{1, \dots, K\}$, let $h := \tau_n - \tau_{n-1} \le \frac{s_{\max}}{N}$, and define $I := \{i \mid x_i \ge \tau_n\}$ and $J := \{j \mid -y_j \ge \tau_n\}$. Since neither $\max(x_i, 0)$ nor $\min(-y_j, s_{\max})$ lies in $(\tau_{n-1}, \tau_n)$, we have $x_i(\tau_n) - x_i(\tau_{n-1}) = h$ for $i \in I$ (with $x_i(\tau_n) = \tau_n$) and $0$ for $i \notin I$. Similarly, $y_j(\tau_n) - y_j(\tau_{n-1}) = -h$ for $j \in J$ (with $y_j(\tau_n) = -\tau_n$) and $0$ for $j \notin J$ (with $y_j(\tau_n) = y_j(\tau_{n-1}) = y_j \ge -\tau_{n-1}$).

\begin{claim}
\label{clm:Lstep}
$L(\tau_n) - L(\tau_{n-1}) \le h\sum_{i \in I}\sum_{j \notin J} W_{ij}(\tau_n + y_j)^2$.
\end{claim}

\begin{proof}
Since $C(\tau_n) - C(\tau_{n-1}) = h\sum_{i \in I} u_i u_i^* - h\sum_{j \in J} v_j v_j^*$ and $h > 0$, applying Lemma~\ref{lem:subgrad} to $X=C(\tau_{n-1}), Y= C(\tau_n)$ and using $\tr(C(\tau_n)_+^2 z z^*) = \|C(\tau_n)_+ z\|_2^2$ gives
\begin{align*}
L(\tau_n) - L(\tau_{n-1})
&\le \tr\Bigl(C(\tau_n)_+^2\bigl(C(\tau_n) - C(\tau_{n-1})\bigr)\Bigr) \\
&= h\sum_{i \in I}\|C(\tau_n)_+ u_i\|_2^2 - h\sum_{j \in J}\|C(\tau_n)_+ v_j\|_2^2
\le h\sum_{i \in I}\|C(\tau_n)_+ u_i\|_2^2.
\end{align*}
For each $i \in I$, since $x_i(\tau_n) = \tau_n$ and $u_i = \sum_{j=1}^d (v_j^* u_i) v_j$, we have
\[
C(\tau_n) u_i
= \tau_n u_i + \sum_{j=1}^d y_j(\tau_n)(v_j^* u_i) v_j
= \sum_{j=1}^d \bigl(\tau_n + y_j(\tau_n)\bigr)(v_j^* u_i) v_j
= \sum_{j \notin J} (\tau_n + y_j)(v_j^* u_i) v_j,
\]
where the last step uses $y_j(\tau_n) = -\tau_n$ for $j \in J$ and $y_j(\tau_n) = y_j$ for $j \notin J$. Since $(v_1, \dots, v_d)$ is orthonormal and $|v_j^* u_i|^2 = W_{ij}$, using $\|C(\tau_n)_+ u_i\|_2 \le \|C(\tau_n) u_i\|_2$ gives
\[
\|C(\tau_n)_+ u_i\|_2^2 \le \|C(\tau_n) u_i\|_2^2 = \sum_{j \notin J} (\tau_n + y_j)^2 |v_j^* u_i|^2 = \sum_{j \notin J} W_{ij}(\tau_n + y_j)^2.
\]
Summing over $i \in I$ and multiplying by $h > 0$ yields the claim.
\end{proof}
\begin{claim}
\label{clm:Rstep}
$R(\tau_n) - R(\tau_{n-1}) \ge h\sum_{i \in I}\sum_{j \notin J} W_{ij}(\tau_n + y_j)^2 - d s_{\max} h^2$.
\end{claim}
\begin{proof}
We analyze $(x_i(\tau_n) + y_j(\tau_n))_+^3 - (x_i(\tau_{n-1}) + y_j(\tau_{n-1}))_+^3$ for each $(i, j)$ in three cases:

(i) $j \in J$: Then $y_j \le -\tau_n < -\tau_{n-1}$, so $y_j(\tau_n) = -\tau_n$ and $y_j(\tau_{n-1}) = -\tau_{n-1}$. Since $x_i(s) \le s$ for all $s \in [0, s_{\max}]$, we have $x_i(\tau_n) + y_j(\tau_n) \le 0$ and $x_i(\tau_{n-1}) + y_j(\tau_{n-1}) \le 0$, and hence $(x_i(\tau_n) + y_j(\tau_n))_+ = (x_i(\tau_{n-1}) + y_j(\tau_{n-1}))_+ = 0$.

(ii) $j \notin J, i \notin I$: Then $x_i(\tau_n) = x_i(\tau_{n-1})$ and $y_j(\tau_n) = y_j(\tau_{n-1})$, so the two terms are equal and their difference is $0$.

(iii) $j \notin J, i \in I$: Then $x_i(\tau_n) = \tau_n$, $x_i(\tau_{n-1}) = \tau_{n-1} = \tau_n - h$, and $y_j(\tau_n) = y_j(\tau_{n-1}) = y_j \ge -\tau_{n-1}$. Hence $\tau_n + y_j - h = \tau_{n-1} + y_j \ge 0$ and $0 < \tau_n + y_j \le \tau_n \le s_{\max}$ (since $h > 0$ and $y_j \le 0$), which gives
\begin{align*}
\frac{1}{3}\Bigl(\bigl(x_i(\tau_n) + y_j(\tau_n)\bigr)_+^3 - \bigl(x_i(\tau_{n-1}) + y_j(\tau_{n-1})\bigr)_+^3\Bigr)
&= \frac{1}{3}\bigl((\tau_n + y_j)^3 - (\tau_n + y_j - h)^3\bigr) \\
&= h(\tau_n + y_j)^2 - h^2(\tau_n + y_j) + \frac{1}{3}h^3 \\
&\ge h(\tau_n + y_j)^2 - s_{\max} h^2.
\end{align*}

Multiplying by $W_{ij}$ and summing over $i, j$, we obtain
\begin{align*}
R(\tau_n) - R(\tau_{n-1})
&= \frac{1}{3}\sum_{i,j=1}^d W_{ij}\Bigl(\bigl(x_i(\tau_n) + y_j(\tau_n)\bigr)_+^3 - \bigl(x_i(\tau_{n-1}) + y_j(\tau_{n-1})\bigr)_+^3\Bigr) \\
&= \frac{1}{3}\sum_{i \in I}\sum_{j \notin J} W_{ij}\bigl((\tau_n + y_j)^3 - (\tau_n + y_j - h)^3\bigr) \\
&\ge h\sum_{i \in I}\sum_{j \notin J} W_{ij}(\tau_n + y_j)^2 - s_{\max} h^2 \sum_{i \in I}\sum_{j \notin J} W_{ij} \\
&\ge h\sum_{i \in I}\sum_{j \notin J} W_{ij}(\tau_n + y_j)^2 - s_{\max} h^2 \sum_{i,j=1}^d W_{ij} \\
&= h\sum_{i \in I}\sum_{j \notin J} W_{ij}(\tau_n + y_j)^2 - d s_{\max} h^2,
\end{align*}
where the last equality uses $\sum_{j=1}^d W_{ij} = \sum_{j=1}^d |u_i^* v_j|^2 = 1$ for each $i$.
\end{proof}

\begin{proof}[Proof of \Cref{lm:maintech}]
If $s_{\max} = 0$, then \Cref{clm:endpoint} gives $\frac{1}{3}\tr(C_+^3) \le L(0) = 0 = R(0) = R(s_{\max}) = \frac{1}{3}\sum_{i,j=1}^d W_{ij}(\lambda_i + \mu_j)_+^3$, which proves~\eqref{eq:peano} since $C = A + B$.

Suppose $s_{\max} > 0$. For any integer $N \ge 1$ and each $n \in \{1, \dots, K\}$, combining \Cref{clm:Lstep} and \Cref{clm:Rstep} (with $h = \tau_n - \tau_{n-1}$) yields
\[
L(\tau_n) - L(\tau_{n-1}) \le R(\tau_n) - R(\tau_{n-1}) + d s_{\max}(\tau_n - \tau_{n-1})^2.
\]
Summing over $n = 1, \dots, K$ and using $L(0) = R(0) = 0$ and $\tau_n - \tau_{n-1} \le \frac{s_{\max}}{N}$, we obtain
\begin{align*}
L(s_{\max})
= \sum_{n=1}^K \bigl(L(\tau_n) - L(\tau_{n-1})\bigr)
&\le \sum_{n=1}^K \bigl(R(\tau_n) - R(\tau_{n-1})\bigr) + d s_{\max}\sum_{n=1}^K (\tau_n - \tau_{n-1})^2 \\
&\le R(s_{\max}) + \frac{d s_{\max}^2}{N}\sum_{n=1}^K (\tau_n - \tau_{n-1})
= R(s_{\max}) + \frac{d s_{\max}^3}{N}.
\end{align*}
Since $N \ge 1$ is arbitrary, letting $N \to \infty$ gives $L(s_{\max}) \le R(s_{\max})$, and \Cref{clm:endpoint} then yields
\[
\frac{1}{3}\tr(C_+^3) \le L(s_{\max}) \le R(s_{\max}) = \frac{1}{3}\sum_{i,j=1}^d W_{ij}(\lambda_i + \mu_j)_+^3.
\]
Recalling $C = A + B$ proves~\eqref{eq:peano}.
\end{proof}

\section{Extension}
\label{sec:extension}

Replacing the monomial integral in the proof of \Cref{thm:binomial} by a piecewise-linear interpolation of $f''$ extends \Cref{lm:maintech} to all functions with convex second derivative.

\begin{theorem}[Fourth-order trace inequality]
\label{thm:general}
Let $A, B \in \Herm^d$ have eigendecompositions $A = \sum_{i=1}^d \lambda_i u_i u_i^*$ and $B = \sum_{j=1}^d \mu_j v_j v_j^*$, let $W_{ij} := |u_i^* v_j|^2$, and let
\[
a := \min_{1 \le i, j \le d}(\lambda_i + \mu_j), \qquad b := \max_{1 \le i, j \le d}(\lambda_i + \mu_j).
\]
Let $f : [a, b] \to \R$ be continuous on $[a, b]$ and twice continuously differentiable on $(a, b)$. If $f''$ is convex on $(a, b)$, then
\begin{equation}
\label{eq:general}
\tr(f(A + B)) \le \sum_{i,j=1}^d W_{ij} f(\lambda_i + \mu_j).
\end{equation}
If $f''$ is concave on $(a, b)$, then the reverse inequality holds.
\end{theorem}

Here $f(A + B)$ is well defined because every eigenvalue of $A + B$ lies in $[a, b]$: $a$ and $b$ are the sums of the smallest, respectively largest, eigenvalues of $A$ and $B$. The convexity hypothesis holds in particular when $f$ is three times differentiable on $(a, b)$ with $f'''$ nondecreasing (then $f''$ is continuous and convex), and hence whenever $f^{(4)} \ge 0$ on $(a, b)$.

\ignore{Replacing the monomial integral in the proof of \Cref{thm:binomial} by piecewise-linear spline interpolation of $f''$ extends \Cref{lm:maintech} to all functions with convex second derivative. Throughout this section, we use the same notation as in \Cref{sec:main}: $A = \sum_{i=1}^d \lambda_i u_i u_i^* \in \Herm^d$, $B = \sum_{j=1}^d \mu_j v_j v_j^* \in \Herm^d$, $W_{ij} := |u_i^* v_j|^2$, and $\sigma_1, \dots, \sigma_d$ are the eigenvalues of $A + B$.
Let $a := \min_{1 \leq i,j \leq d} (\lambda_i + \mu_j)$
and $b := \max_{1 \leq i,j \leq d} (\lambda_i + \mu_j)$.  Note
that by Rayleigh quotients, $\sigma_i \in [a, b]$
for all $i$.

\begin{theorem}[Fourth-order trace inequality]
\label{thm:general}
Let $f \in C([a, b])$. If $f \in C^2((a, b))$ with $f''$ convex on $(a, b)$ (in particular, if $f$ is three times differentiable on $(a, b)$ with $f'''$ nondecreasing on $(a, b)$), then
\begin{equation}
\label{eq:general}
\tr(f(A + B)) \le \sum_{i,j=1}^d W_{ij} f(\lambda_i + \mu_j).
\end{equation}
If instead $f \in C^2((a, b))$ with $f''$ concave on $(a, b)$ (in particular, if $f$ is three times differentiable on $(a, b)$ with $f'''$ nonincreasing on $(a, b)$), then the reverse inequality holds in~\eqref{eq:general}.
\end{theorem}}
\begin{proof}
Let $\sigma_1, \dots, \sigma_d$ be the eigenvalues of $A + B$; as noted above, $\sigma_i \in [a, b]$. By \Cref{lem:overlap}(ii) and linearity, equality holds in~\eqref{eq:general} for every cubic polynomial. If $a = b$, then $\sigma_i = \lambda_i + \mu_j = a$ for all $i, j$, so both sides of~\eqref{eq:general} equal $d f(a)$ since $\sum_{i,j=1}^d W_{ij} = d$. Suppose now that $a < b$. % Note first that if $f$ is three times differentiable on $(a, b)$ with $f'''$ nondecreasing, then $f''$ is differentiable (hence continuous) on $(a, b)$ with nondecreasing derivative, so $f \in C^2((a, b))$ and $f''$ is convex on $(a, b)$.

Fix $\rho \in (0, 1)$, and for $x \in [a, b]$ define the rescaled function $f_\rho(x) := f\bigl(\rho x + (1 - \rho)\frac{a + b}{2}\bigr)$. Since $\rho x + (1 - \rho)\frac{a + b}{2} \in (a, b)$ for all $x \in [a, b]$, we have $f_\rho \in C^2([a, b])$, and $\psi := f_\rho'' \in C([a, b])$ is convex on $[a, b]$. Integrating by parts twice gives, for all $x \in [a, b]$,
\[
f_\rho(x) = f_\rho(a) + f_\rho'(a)(x - a) + \int_a^x (x - s)\psi(s)\, ds.
\]

For each integer $N \ge 2$, let $h := \frac{b - a}{N}$ and $\tau_n := a + n h$ for $n \in \{0, \dots, N\}$, and define
\begin{align*}
c_\ell &:= \frac{\psi(\tau_{\ell+1}) - 2\psi(\tau_\ell) + \psi(\tau_{\ell-1})}{h} && (1 \le \ell \le N - 1), \\
\psi_N(s) &:= \psi(\tau_0) + \frac{\psi(\tau_1) - \psi(\tau_0)}{h}(s - a) + \sum_{\ell=1}^{N-1} c_\ell(s - \tau_\ell)_+ && (s \in [a, b]).
\end{align*}
Since $\tau_\ell = \frac{\tau_{\ell-1} + \tau_{\ell+1}}{2}$ and $\psi$ is convex on $[a, b]$, we have $c_\ell \ge 0$ for all $\ell$. On each subinterval $[\tau_{n-1}, \tau_n]$, for $n \in \{1, \ldots, N\}$, $\psi_N$ is affine with slope $\frac{\psi(\tau_1) - \psi(\tau_0)}{h} + \sum_{\ell=1}^{n-1} c_\ell = \frac{\psi(\tau_n) - \psi(\tau_{n-1})}{h}$; since $\psi_N(\tau_0) = \psi(\tau_0)$, it follows that $\psi_N(\tau_n) = \psi(\tau_n)$ for all $n \in \{0, \dots, N\}$. In particular, for $s \in [\tau_{n-1}, \tau_n]$, $\psi_N(s)$ lies between $\psi(\tau_{n-1})$ and $\psi(\tau_n)$, so $|\psi(s) - \psi_N(s)| \le \max\bigl(|\psi(s) - \psi(\tau_{n-1})|, |\psi(s) - \psi(\tau_n)|\bigr)$. Since $\tau_n - \tau_{n-1} = h = \frac{b - a}{N} \to 0$ and $\psi$ is uniformly continuous on $[a, b]$, we obtain $\max_{s \in [a, b]} |\psi(s) - \psi_N(s)| \to 0$ as $N \to \infty$.

Next, for $x \in [a, b]$, define 
\[
f_{\rho, N}(x) := f_\rho(a) + f_\rho'(a)(x - a) + \int_a^x (x - s)\psi_N(s)\, ds.
\]
Since $\tau_\ell > a$ and $\int_a^x (x - s)(s - \tau_\ell)_+\, ds = \frac{1}{6}(x - \tau_\ell)_+^3$ for all $x \in [a, b]$, integrating $\psi_N$ term by term gives
\[
f_{\rho, N}(x) = P_N(x) + \frac{1}{6}\sum_{\ell=1}^{N-1} c_\ell(x - \tau_\ell)_+^3,
\]
where $P_N(x) := f_\rho(a) + f_\rho'(a)(x - a) + \frac{1}{2}\psi(\tau_0)(x - a)^2 + \frac{\psi(\tau_1) - \psi(\tau_0)}{6 h}(x - a)^3$ is a cubic polynomial. Comparing the integral representations of $f_\rho$ and $f_{\rho, N}$ gives, for all $x \in [a, b]$,
\[
|f_\rho(x) - f_{\rho, N}(x)| = \left|\int_a^x (x - s)\bigl(\psi(s) - \psi_N(s)\bigr)\, ds\right| \le \frac{(b - a)^2}{2}\max_{s \in [a, b]} |\psi(s) - \psi_N(s)|,
\]
so $\max_{x \in [a, b]} |f_\rho(x) - f_{\rho, N}(x)| \to 0$ as $N \to \infty$.

Since $\sigma_i \in [a, b]$ and $\lambda_i + \mu_j \in [a, b]$ for all $i, j$, applying \Cref{lem:overlap}(ii) to $P_N$ and \Cref{lm:maintech} to $A, B - \tau_\ell I_d$ for each $\ell \in \{1, \dots, N - 1\}$ (using $c_\ell \ge 0$) yields
\[
\tr(f_{\rho, N}(A + B)) \le \sum_{i,j=1}^d W_{ij} f_{\rho, N}(\lambda_i + \mu_j).
\]
Letting $N \to \infty$ gives $\tr(f_\rho(A + B)) \le \sum_{i,j=1}^d W_{ij} f_\rho(\lambda_i + \mu_j)$, and then letting $\rho \to 1$ from below  (using continuity of $f$ on $[a, b]$) proves~\eqref{eq:general}. Finally, applying~\eqref{eq:general} to $-f$ yields the reverse inequality when $f''$ is concave.
\end{proof}

\section{Applications}
\label{sec:applications}
We now present four applications of \Cref{thm:general} in which the overlaps $W_{ij}$ disappear from the final statements. We keep the notation of \Cref{sec:main}. For $M \in \Herm^d$ with eigenvalues $\nu_1, \dots, \nu_d$, we write $|M| := (M^2)^{1/2}$, $\|M\|_{S_p} := \tr(|M|^p)^{1/p} = (\sum_{i=1}^d |\nu_i|^p)^{1/p}$ ($p \ge 1$) for its Schatten-$p$ norm, and $\|M\|_{\mathrm{op}} := \max_{1 \le i \le d} |\nu_i|$ for its operator norm.

\subsection{Geometry of Schatten-\texorpdfstring{$p$}{p} Spaces: Two-Point and Hanner Inequalities}
\label{sec:schatten}
For $p \ge 3$, \Cref{thm:general} applies to $f(x) = |x|^p$ (even though $f \notin C^4(\R)$ for $3 \le p < 4$) and gives $$\|A \pm B\|_{S_p}^p \le \sum_{i,j} W_{ij}|\lambda_i \pm \mu_j|^p.$$ 

Now, for $q \ge 1$ and a real array $x = (x_{ij})_{i,j=1}^d$, write $\|x\|_{\ell_q(W)} := \bigl(\sum_{i,j=1}^d W_{ij}|x_{ij}|^q\bigr)^{1/q}$ for the $\ell_q$ norm on $\{1, \dots, d\}^2$ weighted by the overlaps. 

Since $\sum_{j=1}^d W_{ij} = \sum_{i=1}^d W_{ij} = 1$, the arrays $x_{ij} := \lambda_i$ and $y_{ij} := \mu_j$ satisfy $\|x\|_{\ell_q(W)} = \|A\|_{S_q}$ and $\|y\|_{\ell_q(W)} = \|B\|_{S_q}$. Thus, any upper bound on $\|x + y\|^p + \|x - y\|^p$ valid in $\ell_p$ spaces transfers to $S_p$. This gives short proofs, for $p \ge 3$, of the two-point inequality of~\cite{ball-carlen-lieb} and of Hanner's \cite{hanner1956} inequality in $S_p$.

%For $3 \le p < 4$, the function $f(x) = |x|^p$ is not in $C^4(\R)$, yet \Cref{thm:general} applies across all $p \ge 3$.  Since the overlap weights $W_{ij}$ preserve individual Schatten-$p$ norms while dominating $\|A \pm B\|_{S_p}$, they immediately transfer $L_p$-smoothness and uniform convexity inequalities over vectors to $S_p$, giving short proofs of both the sharp two-point ($2$-smoothness) inequality of~\cite{ball-carlen-lieb} and the inequality of~\cite{hanner1956} in $S_p$ for all $p \ge 3$.

\begin{corollary}[Schatten-$p$ two-point inequality]
\label{cor:twopoint}
For every $p \ge 3$ and all $A, B \in \Herm^d$,
\[
\left(\frac{\|A + B\|_{S_p}^p + \|A - B\|_{S_p}^p}{2}\right)^{\!2/p} \le \|A\|_{S_p}^2 + (p - 1)\|B\|_{S_p}^2.
\]
\end{corollary}

\begin{proof}
For every $p \ge 3$, $f(x) = |x|^p$ is in $C^2(\R)$ with convex second derivative $f''(x) = p(p - 1)|x|^{p-2}$.  \Cref{thm:general} therefore yields
\begin{equation}
\label{eq:sp-dominance}
\|A + B\|_{S_p}^p = \Tr\bigl(|A + B|^p\bigr) \le \sum_{i,j=1}^d W_{ij}|\lambda_i + \mu_j|^p.
\end{equation}

By the Bonami--Gross scalar two-point inequality~\citep{bonami, gross} (see also~\citep{ball-carlen-lieb}),
\[
\frac{|x + y|^p + |x - y|^p}{2} \le \bigl(x^2 + (p - 1)y^2\bigr)^{p/2} \qquad (x, y \in \R,\; p \ge 2).
\]
Applying~\eqref{eq:sp-dominance} to $(A, B)$ and to $(A, -B)$ (which share the eigenbasis overlaps $W_{ij}$) and using this inequality yields
\[
\frac{\|A + B\|_{S_p}^p + \|A - B\|_{S_p}^p}{2} \le \sum_{i,j=1}^d W_{ij}\frac{|\lambda_i + \mu_j|^p + |\lambda_i - \mu_j|^p}{2} \le \sum_{i,j=1}^d W_{ij}\bigl(\lambda_i^2 + (p - 1)\mu_j^2\bigr)^{p/2}.
\]
Since $\frac{p}{2} \ge \frac{3}{2} > 1$ and $\sum_{j=1}^d W_{ij} = \sum_{i=1}^d W_{ij} = 1$, Minkowski's triangle inequality in $\ell_{p/2}(W)$ gives
\begin{align*}
\left(\sum_{i,j=1}^d W_{ij}\bigl(\lambda_i^2 + (p - 1)\mu_j^2\bigr)^{p/2}\right)^{\!2/p}
&\le \left(\sum_{i,j=1}^d W_{ij}|\lambda_i|^p\right)^{\!2/p} + (p - 1)\left(\sum_{i,j=1}^d W_{ij}|\mu_j|^p\right)^{\!2/p} \\
&= \|A\|_{S_p}^2 + (p - 1)\|B\|_{S_p}^2.\qedhere
\end{align*}
\end{proof}

\begin{corollary}[Hanner's inequality in Schatten-$p$]
\label{cor:hanner}
For every $p \ge 3$ and all $A, B \in \Herm^d$,
\[
\|A + B\|_{S_p}^p + \|A - B\|_{S_p}^p \le \bigl(\|A\|_{S_p} + \|B\|_{S_p}\bigr)^p + \bigl|\|A\|_{S_p} - \|B\|_{S_p}\bigr|^p.
\]
\end{corollary}
\begin{proof}
Applying~\eqref{eq:sp-dominance} to $(A, B)$ and $(A, -B)$ gives
\[
\|A + B\|_{S_p}^p + \|A - B\|_{S_p}^p \le \sum_{i,j=1}^d W_{ij}\bigl(|\lambda_i + \mu_j|^p + |\lambda_i - \mu_j|^p\bigr).
\]
In $\ell_p(W)$ on $\{1, \dots, d\}^2$, the arrays $x_{ij} := \lambda_i$ and $y_{ij} := \mu_j$ satisfy $\|x\|_{\ell_p(W)} = \bigl(\sum_{i,j=1}^d W_{ij}|\lambda_i|^p\bigr)^{1/p} = \|A\|_{S_p}$ and $\|y\|_{\ell_p(W)} = \|B\|_{S_p}$, while the RHS above equals $\|x + y\|_{\ell_p(W)}^p + \|x - y\|_{\ell_p(W)}^p$. Applying Hanner's classical inequality in $\ell_p(W)$~\citep{hanner1956} completes the proof.
\end{proof}

\citet{ball-carlen-lieb} proved Hanner's inequality in $S_p$ for $p \ge 4$ (and by duality for $1 \le p \le \frac{4}{3}$); \cite{heinavaara2024} resolved $p = 3$ via Koldobsky's Fourier characterization of $L_3$-embeddability and later all $p \ge 1$~\citep{heinavaara2025} via tracial joint spectral measures. \Cref{cor:hanner} gives an elementary proof across all $p \ge 3$ (covering $3 \le p < 4$), and standard block Hermitization and H\"older duality~\cite[Lemma~6]{ball-carlen-lieb} extend \Cref{cor:twopoint,cor:hanner} to general complex matrices and yield the corresponding reversed inequalities in the dual range $1 \le p \le \frac{3}{2}$.

\subsection{Strengthening Clarkson--McCarthy Convexity on the PSD Cone}
\label{sec:clarkson}

Recall the classical Clarkson--McCarthy convexity inequality~\citep{dixmier1953, mccarthy1967}: for all $A, B \in \Herm^d$ and $p \ge 2$,
\begin{equation}
\label{eq:clarkson-classic}
\frac{\|A\|_{S_p}^p + \|B\|_{S_p}^p}{2} - \left\|\frac{A + B}{2}\right\|_{S_p}^p \ge \left\|\frac{A - B}{2}\right\|_{S_p}^p.
\end{equation}
On all of $\Herm^d$, the inequality~\eqref{eq:clarkson-classic} is tight (with equality at $B = -A$). Restricting to PSD matrices $A, B \in \Herm_+^d$, \Cref{thm:general} yields a stronger Frobenius-norm bound.

\begin{corollary}[PSD Clarkson--McCarthy inequality]
\label{cor:strong-convex}
For every $p \ge 3$ and all $A, B \in \Herm_+^d$,
\begin{equation}
\label{eq:hs-strong}
\frac{\Tr(A^p) + \Tr(B^p)}{2} - \Tr\!\left(\left(\frac{A + B}{2}\right)^{\!p}\right)
\ge \frac{1}{2}\!\left(1 - \frac{1}{p} \right)\,\bigl\|A^{p/2} - B^{p/2}\bigr\|_F^2.
\end{equation}
\end{corollary}
\begin{proof}
We first show the scalar inequality
\begin{equation}
\label{eq:scalar-cm}
\left(\frac{x + y}{2}\right)^p \le \frac{1}{p}\frac{x^p + y^p}{2} + \left(1 - \frac{1}{p}\right)(xy)^{p/2} \qquad (x, y \ge 0,\; p \ge 2).
\end{equation}
If $x + y = 0$, both sides are $0$. Otherwise, wlog assume $x \ge y \ge 0$, let $s := \frac{x + y}{2} > 0$ and $t := \frac{x - y}{x + y} \in [0, 1]$ (so $x = s(1 + t)$ and $y = s(1 - t)$), and divide~\eqref{eq:scalar-cm} by $s^p$ to rewrite it as $F(t) \ge 1$, where $F(t) := \frac{(1 + t)^p + (1 - t)^p}{2p} + \bigl(1 - \frac{1}{p}\bigr)(1 - t^2)^{p/2}$. Note that $F \in C([0, 1])$ with $F(0) = 1$. For $t \in (0, 1)$, setting $e^{2w} := \frac{1 + t}{1 - t}$ ($w > 0$, so $1 \pm t = \sqrt{1 - t^2}\,e^{\pm w}$ and $t = \sqrt{1 - t^2}\sinh w$) gives
\[
F'(t) = \frac{(1 + t)^{p - 1} - (1 - t)^{p - 1}}{2} - (p - 1)t(1 - t^2)^{\frac{p - 2}{2}} = (1 - t^2)^{\frac{p - 1}{2}}\Bigl(\sinh\bigl((p - 1)w\bigr) - (p - 1)\sinh w\Bigr) \ge 0,
\]
since $\sinh(rw) = \sum_{k=0}^\infty \frac{(rw)^{2k+1}}{(2k+1)!} \ge r\sinh w$ for $r := p - 1 \ge 1$ and $w > 0$. Thus $F(t) \ge F(0) = 1$ on $[0, 1]$, proving~\eqref{eq:scalar-cm}.

Applying \Cref{thm:general} with $f(x) = (x/2)^p$ ($p \ge 3$) on $[0, \infty)$ (where $f''(x) = \frac{p(p-1)}{4}(x/2)^{p-2}$ is convex) and summing~\eqref{eq:scalar-cm} against $W_{ij}$ at $(\lambda_i, \mu_j)$ gives
\[
\tr\!\left(\left(\frac{A + B}{2}\right)^p\right) \le \frac{1}{p}\frac{\tr(A^p) + \tr(B^p)}{2} + \left(1 - \frac{1}{p}\right)\tr(A^{p/2} B^{p/2}),
\]
since $\sum_{i,j=1}^d W_{ij}\lambda_i^{p/2}\mu_j^{p/2} = \tr(A^{p/2} B^{p/2})$. Substituting $\tr(A^{p/2} B^{p/2}) = \frac{1}{2}\bigl(\tr(A^p) + \tr(B^p) - \|A^{p/2} - B^{p/2}\|_F^2\bigr)$ 
completes the proof.
\end{proof}
Since $x \mapsto x^{2/p}$ is operator monotone on $[0, \infty)$ for $p \ge 2$, the norm inequality of~\citet{bks1975} and~\citet{ando1988} gives $\|U^{2/p} - V^{2/p}\|_{S_p} \le \||U - V|^{2/p}\|_{S_p} = \|U - V\|_F^{2/p}$ for all $U, V \in \Herm_+^d$. Taking $U = A^{p/2}$ and $V = B^{p/2}$ yields $\|A - B\|_{S_p}^p \le \|A^{p/2} - B^{p/2}\|_F^2$, so~\eqref{eq:hs-strong} implies
\[
\frac{\Tr(A^p) + \Tr(B^p)}{2} - \Tr\!\left(\left(\frac{A + B}{2}\right)^{\!p}\right)
\ge \frac{(p - 1)2^{p-1}}{p}\left\|\frac{A - B}{2}\right\|_{S_p}^p,
\]
strengthening~\eqref{eq:clarkson-classic} on $\Herm_+^d$ everywhere by the factor $\frac{(p - 1)2^{p-1}}{p}$. Two regimes illustrate the gain:
\begin{itemize}[nosep]
\item \emph{Near the identiy:} For $X, Y \in \Herm^d$ and small $t \in \R$, set $A = I_d + tX$ and $B = I_d + tY$. The right-hand side of~\eqref{eq:clarkson-classic} is only $O(|t|^p)$, whereas both sides of~\eqref{eq:hs-strong} equal $\frac{p(p - 1)}{8}t^2\|X - Y\|_F^2 + O(|t|^3)$. In particular, the constant $\frac{1}{2}(1 - \frac{1}{p})$ cannot be improved.
%\item \emph{Near the identity:} For $X, Y \in \Herm^d$ and sufficiently small $t \in \R$, setting $A = I_d + tX$ and $B = I_d + tY$ in~\eqref{eq:clarkson-classic} gives only $\frac{|t|^p}{2^p}\|X - Y\|_{S_p}^p = O(|t|^p)$, whereas both sides of~\eqref{eq:hs-strong} expand to $\frac{p(p - 1)}{8}t^2\|X - Y\|_F^2 + \frac{p(p - 1)(p - 2)}{16}t^3\Tr((X - Y)(X^2 - Y^2)) + O(t^4)$ (matching through $O(t^3)$ and showing that the constant $\frac{1}{2}(1 - \frac{1}{p})$ is sharp).
\item \emph{Orthogonal projections:} For $d \geq 2$, if $A, B \in \Herm_+^d$ are orthogonal rank-$1$ projections ($A^2 = A$, $B^2 = B$, $AB = 0$), then~\eqref{eq:clarkson-classic} gives only $\Tr((\frac{A+B}{2})^p) \le 1 - 2^{1-p} \to 1$ as $p \to \infty$, whereas~\eqref{eq:hs-strong} gives $\Tr((\frac{A+B}{2})^p) \le \frac{1}{p} \to 0$.
\end{itemize}

\subsection{Hyperbolic and Trigonometric Trace Inequalities}
\label{sec:trig}

\Cref{thm:general} also yields matrix trace analogues of the classical addition formulas for $\cosh$, $\sinh$, $\cos$, and $\sin$, together with unconditional commutator bounds for $\cos$ and $\sin$.

\begin{corollary}[Hyperbolic and trigonometric trace inequalities]
\label{cor:trig}
Let $A, B \in \Herm^d$.
\begin{enumerate}[label=(\roman*),nosep]
\item 
$\Tr\bigl(\cosh(A + B)\bigr) \le \Tr(\cosh A\cosh B + \sinh A\sinh B).$
\item If $\|A\|_{\mathrm{op}} + \|B\|_{\mathrm{op}} \le \frac{\pi}{2}$, then
$\Tr\bigl(\cos(A + B)\bigr) \le \Tr(\cos A\cos B - \sin A\sin B).$
\item If $A, B \in \Herm_+^d$, then
$\Tr\bigl(\sinh(A + B)\bigr) \le \Tr(\sinh A\cosh B + \cosh A\sinh B).$
\item If $A, B \in \Herm_+^d$ and $\|A\|_{\mathrm{op}} + \|B\|_{\mathrm{op}} \le \pi$, then
$\Tr\bigl(\sin(A + B)\bigr) \le \Tr(\sin A\cos B + \cos A\sin B).$
\item 
$\Bigl|\Tr\bigl(\cos(A + B)\bigr) - \Tr(\cos A\cos B - \sin A\sin B)\Bigr| \le \frac{1}{24}\,\bigl\|[A, B]\bigr\|_F^2,$
and \\ $\Bigl|\Tr\bigl(\sin(A + B)\bigr) - \Tr(\sin A\cos B + \cos A\sin B)\Bigr| \le \frac{1}{24}\,\bigl\|[A, B]\bigr\|_F^2.$
\end{enumerate}
\end{corollary}

\begin{proof}
Expanding in the eigenbases of $A$ and $B$ and using the scalar addition formulas gives $\sum_{i,j=1}^d W_{ij} f(\lambda_i + \mu_j)$ equal to the RHS in each of (i)--(iv); for instance,
\[
\sum_{i,j=1}^d W_{ij}\cos(\lambda_i + \mu_j) = \sum_{i,j=1}^d W_{ij}(\cos\lambda_i\cos\mu_j - \sin\lambda_i\sin\mu_j) = \Tr(\cos A\cos B - \sin A\sin B),
\]
and identically for $\cosh$, $\sinh$, and $\sin$.  Note that the interval $[a, b]$ in \Cref{thm:general} lies in $\R$ for (i), $[-\frac{\pi}{2}, \frac{\pi}{2}]$ for (ii), $[0, \infty)$ for (iii), and $[0, \pi]$ for (iv). In each case $f^{(4)} = f \ge 0$ on $[a, b]$ (so $f''$ is convex), and \Cref{thm:general} yields (i)–(iv). 

For~(v), the functions $g_\pm(x) := \frac{x^4}{24} \pm \cos x$ and $h_\pm(x) := \frac{x^4}{24} \pm \sin x$ satisfy $g_\pm^{(4)}(x) = 1 \pm \cos x \ge 0$ and $h_\pm^{(4)}(x) = 1 \pm \sin x \ge 0$ for all $x \in \R$. Applying \Cref{thm:general} to $g_\pm$ and $h_\pm$ bounds each
LHS in~(v) by
\begin{align*}
\frac{1}{24}\left(\sum_{i,j=1}^d W_{ij}(\lambda_i + \mu_j)^4 - \Tr\bigl((A + B)^4\bigr)\right)
= \frac{1}{12}\Tr(A^2 B^2 - ABAB) = -\frac{1}{24}\Tr\bigl([A, B]^2\bigr) = \frac{1}{24}\,\bigl\|[A, B]\bigr\|_F^2,
\end{align*}
where we used \Cref{lem:overlap}(i) and $[A, B]^* = -[A, B]$.
\end{proof}

While (i) also follows by averaging the Golden--Thompson inequality (\Cref{cor:gt}) at $(A, B)$ and $(-A, -B)$, inequalities (ii)--(v) do not seem to follow from Golden--Thompson. The norm thresholds $\frac{\pi}{2}$ in (ii) and $\pi$ in (iv) are sharp: for $A = \frac{c}{2}I_2 + \varepsilon\begin{psmallmatrix}0&1\\1&0\end{psmallmatrix}$ and $B = \frac{c}{2}I_2 + \varepsilon\begin{psmallmatrix}1&0\\0&-1\end{psmallmatrix}$ ($c \ge 2\varepsilon > 0$, so $A, B \in \Herm_+^2$ with
$\|A\|_{\mathrm{op}} + \|B\|_{\mathrm{op}} = c + 2 \varepsilon$), Taylor expansion as $\varepsilon \to 0$ gives $\sum_{i,j=1}^2 W_{ij} f(\lambda_i + \mu_j) - \Tr(f(A + B)) = \frac{1}{3}f^{(4)}(c)\varepsilon^4 + O(\varepsilon^5)$ (since $\|[A, B]\|_F^2 = 8\varepsilon^4$), which becomes negative whenever $\cos c < 0$ in (ii) or $\sin c < 0$ in (iv).

\subsection{Random Unitary Orbits and Free Additive Convolution}
\label{sec:free}

Averaging \Cref{thm:general} over a Haar-random unitary conjugation replaces $W_{ij}$ by its uniform expectation $\frac{1}{d}$, giving a finite-dimensional orbit inequality whose $d \to \infty$ limit recovers the sufficiency direction of the free convolution comparison theorem~\cite[Theorem~1.3]{Heinavaara2026}.

\begin{corollary}[Haar-unitary orbit and free convolution comparison]
\label{cor:free}
\leavevmode
\begin{enumerate}[label=(\roman*),nosep]
\item \emph{(Finite-$d$ unitary orbit)} Let $A, B \in \Herm^d$, $a := \min_{1 \leq i,j \leq d} (\lambda_i + \mu_j)$,  $b := \max_{1 \leq i,j \leq d} (\lambda_i + \mu_j)$, and let $U$ be drawn from the Haar probability measure on the unitary group $\mathrm{U}(d)$. For every $f \in C([a, b]) \cap C^2((a, b))$ with $f''$ convex on $(a, b)$,
\[
\mathbb{E}_U\!\left[\frac{1}{d}\Tr\bigl(f(A + U B U^*)\bigr)\right] \le \frac{1}{d^2}\sum_{i, j=1}^d f(\lambda_i + \mu_j).
\]
\item \emph{(Free vs.\ classical additive convolution~\citep{Heinavaara2026})} Let $\mu$ and $\nu$ be compactly supported Borel probability measures on $\R$ whose supports have convex hulls $[a_1, b_1]$ and $[a_2, b_2]$, and set $[a, b] := [a_1 + a_2, b_1 + b_2]$. For every $f \in C([a, b]) \cap C^2((a, b))$ with $f''$ convex on $(a, b)$,
\[
\int_{[a, b]} f\,d(\mu \boxplus \nu) \le \int_{[a, b]} f\,d(\mu * \nu).
\]
\end{enumerate}
\end{corollary}
\begin{proof}
For (i), applying \Cref{thm:general} to $A$ and $U B U^*$ (whose eigenvectors are $u_i$ and $U v_j$) yields
\[
\frac{1}{d}\tr\bigl(f(A + U B U^*)\bigr) \le \frac{1}{d}\sum_{i,j=1}^d |u_i^* U v_j|^2 f(\lambda_i + \mu_j).
\]
By unitary invariance of Haar measure, $U v_j$ is uniformly distributed on the unit sphere, so
\[
\mathbb{E}_U\bigl[|u_i^* U v_j|^2\bigr] = \frac{1}{d}\sum_{\ell=1}^d \mathbb{E}_U\bigl[|u_\ell^* U v_j|^2\bigr] = \frac{1}{d}\mathbb{E}_U\bigl[\|U v_j\|_2^2\bigr] = \frac{1}{d};
\]
taking expectations proves (i).

For (ii), choose $A_d = \operatorname{diag}(\alpha_{1,d}, \dots, \alpha_{d,d})$ and $B_d = \operatorname{diag}(\beta_{1,d}, \dots, \beta_{d,d})$ with $\alpha_{i,d} \in [a_1, b_1]$ and $\beta_{j,d} \in [a_2, b_2]$ such that the empirical spectral measures $\mu_d := \frac{1}{d}\sum_{i=1}^d \delta_{\alpha_{i,d}}$ and $\nu_d := \frac{1}{d}\sum_{j=1}^d \delta_{\beta_{j,d}}$ converge weakly to $\mu$ and $\nu$ as $d \to \infty$. Let $U_d$ be Haar-distributed on $\mathrm{U}(d)$ and let $\bar{\mu}_d$ be the expected empirical spectral measure of $A_d + U_d B_d U_d^*$. Since $\alpha_{i,d} + \beta_{j,d} \in [a, b]$ for all $i, j$, part (i) applied to $A_d, B_d$ reads
\[
\int_{[a, b]} f\,d\bar{\mu}_d \le \frac{1}{d^2}\sum_{i,j=1}^d f(\alpha_{i,d} + \beta_{j,d}) = \int_{[a, b]} f\,d(\mu_d * \nu_d).
\]
As $d \to \infty$, $\mu_d * \nu_d$ converges weakly to $\mu * \nu$, and by the asymptotic freeness theorem~\citep{voiculescu1991} (see also~\citep[Theorem~5.4.10]{agz2010}), $\bar{\mu}_d$ converges weakly to $\mu \boxplus \nu$. All these measures are supported in $[a, b]$ and $f \in C([a, b])$, so letting $d \to \infty$ yields (ii).

%For (ii), choose $A_d = \operatorname{diag}(\lambda_{1,d}, \dots, \lambda_{d,d})$ and $B_d = \operatorname{diag}(\mu_{1,d}, \dots, \mu_{d,d})$ with $\lambda_{i,d} \in [a_1, b_1]$ and $\mu_{j,d} \in [a_2, b_2]$ whose empirical spectral measures $\mu_d := \frac{1}{d}\sum_{i=1}^d \delta_{\lambda_{i,d}}$ and $\nu_d := \frac{1}{d}\sum_{j=1}^d \delta_{\mu_{j,d}}$ converge weakly to $\mu$ and $\nu$ as $d \to \infty$. Then $\mu_d * \nu_d$ is supported in $[a, b]$ and converges weakly to $\mu * \nu$, while by the  asymptotic freeness theorem~\citep{voiculescu1991} (see also~\citep[Theorem~5.4.10]{agz2010}), the expected empirical spectral measure of $A_d + U_d B_d U_d^*$ ($U_d \sim \operatorname{Haar}(\mathrm{U}(d))$) is supported in $[a, b]$ and converges weakly to $\mu \boxplus \nu$. Since $f \in C([a, b])$, (i) with $d \to \infty$ yields (ii).
\end{proof}

\bibliographystyle{plainnat}
\bibliography{refs}

\end{document}